# An introduction to the concept of function within Descartes's algebra of segments

Nicol Imperi            Enrico Rogora

*Department of Mathematics*
*Università degli studi di Roma La Sapienza*

**Abstract.** *In his Géométrie (1637) Descartes introduces the algebra of segments. This is a fundamental step in the mathematical treatment of variable quantities before the creation of differential calculus. It is an algebra with symbols but without numbers, in which the covariation between geometric variables, constrained by ruler and compass constructions or with other geometric constructions, can be expressed with symbolic equations. By using algebraic manipulations, it is possible to easily deduce the properties of the corresponding geometric constructions, including those that produce graphs of rational functions. We believe that the study of functions through Descartes's algebra can be didactically effective in teaching and learning the concept of function in secondary school. Firstly, it avoids the reference to real numbers; secondly, the interpretation of formulas as geometric constructions and vice versa facilitates the "transition" from functions understood as processes to functions understood as objects.*

## 1. Introduction

In scientific observation of a phenomenon it is often necessary to relate two variable quantities. In ancient natural philosophy the only variable quantities which can be described mathematically are those that vary uniformly, such as the position of stars, while more general variable quantities were described in an exclusively qualitative manner. Starting from the Middle Ages, in the schools of Oxford and Paris mathematical models capable of describing other types of variations began to be developed. The instruments originally used were purely geometric. Before the creation of differential calculus, the most refined mathematical model for treating variable quantities was elaborated by Descartes in the Géométrie. Then the success of differential calculus, developed a few years later by Leibniz and Newton, quickly absorbed the Cartesian point of view, overshadowing some characteristics that, even from a didactic point of view, deserve to be highlighted. Descartes's Géométrie elaborates Apollonius's approach to conics, achieving new and remarkable results, thanks also to the algebraic symbolism that begins to assert itself with Viète. The work of Apollonius, the last of the great Hellenistic mathematical works to be recovered and understood in the West, contains the seeds of many developments that led to a substantial advancement of mathematical knowledge of the following centuries. Among other things, it deals with the study of the symptom of a conic, thus the equation (without symbols) in the geometric algebra (without numbers) of Apollonius[1].

---

[1] The symptom of a parabola is described by Apollonius as follows:

If a cone is cut by a plane through its axis, and also cut by another plane cutting the base of the cone in a straight line perpendicular to the base of the axial triangle, and if further the diameter of the section is parallel to one [lat- eral] side of the axial triangle, and if any straight line is drawn from the section of the cone to its diameter such that this straight line is parallel to the common section of the cutting plane and of the cone's base, then this straight line dropped to the diameter will equal in square to [the rectangular plane] under the straight line from the section's vertex to [the point] where the straight line dropped to the diameter cuts it off and under another straight line which is to the straight line between the angle of the cone and the vertex of the section as the square on the base of the axial triangle to the rectangular plane under the remaining two sides of the triangle. I call such a section a parabola. (Heiberg, 1891, pp. 36–38)

The symptom expresses a particular type of covariation of two segments, constrained through a geometric construction. The constraint, that is the geometric construction expressed in Apollonius's geometric algebra, thus becomes a mathematical object: the symptom. Two curves are identified if they have the same symptom and then the curve itself identifies with the symptom. Following this idea, Apollonius is able to identify conic sections with particular loci that have the same symptom. The different symptoms to which Apollonius reduces the conic sections are described with reference to the different applications of the areas considered by the Pythagorean school, which are divided into parabolic, hyperbolic and elliptical. It is from the relationship with these constructions that Apollonius chooses the names we still use today.

In the algebra of polynomials built on the algebra of numbers. we identify the symptom with an equation. Apollonius does not use numbers, variables or equations to describe the symptom of a curve. In Descartes's algebra of segments, the operations correspond to more general geometric constructions of the applications of the areas considered by Apollonius and can be expressed with the symbolic formalism of Viète. The symptom of a curve becomes, with Descartes, a symbolic equation in the algebra of segments. There are still no numbers and there are no infinite operations, derivatives, integrals and limits, but the algebraic point of view indicated by Descartes will soon include these elements as well.

Descartes' algebra allows us to express the symptom of curves that are graphs of rational functions of x and the square root of x; these graphs are loci of points that can be constructed with a ruler and a compass. By introducing additional tools (Descartes's compass and mechanism) it is possible to build points of the graphs of more general algebraic functions, but we will not deal with this generality in this work which is dedicated to secondary education.

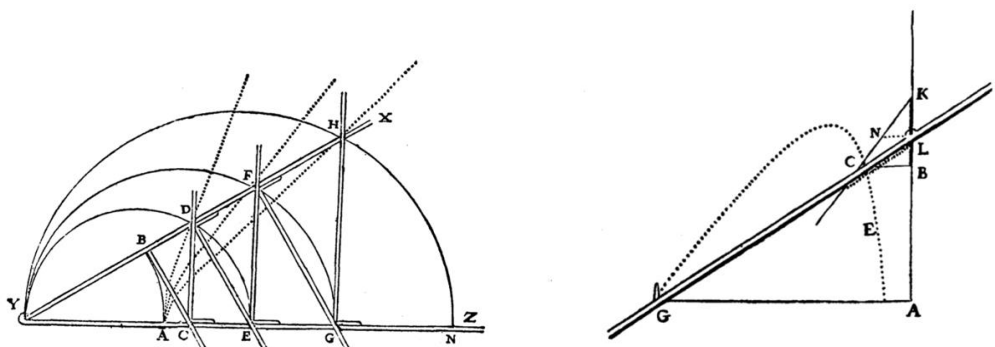

*Figure 1.* These figures describe the tools to extract the roots of any segment (Descartes's Compass on the left, (Descartes, 1637, p. 318)) and to divide an angle into any number of congruent parts (Descartes's mechanism on the right (Descartes, 1637, p. 320)).

We would like to point out the following distinguishing feature of the Cartesian approach to the construction and definition of a curve:
• the modeling of elementary geometric variability with the variable segment (instead of the variable number), which allows a fairly rich and intuitive model without reference to numbers;
• the modeling of the construction process of particular relationships between variable (geometric) quantities through ruler and compass constructions;
• the algebrization of relations, which is accompanied by the use of symbols and symbolic calculation. These are crucial processes in the construction of the concept of function, but in a context that offers a more intuitive and constructive base; a simpler and preparatory context from a didactical point of view, of which we intend to explore the potential.

---

Using algebraic notation and an appropriate reference system, this symptom boils down to the equation $ay = x$.

## 2. The educational potential of Cartesian approach to covariation

In upper secondary schools the definition of a real function of a real variable is introduced as a particular case of the general set-theoretic definition of function. Educational research has identified difficulties related to learning the concept of function as exemplary of some of the most critical aspects in the learning and teaching process of mathematics (Dubinsky & Harel, 1992; Sfard, 1992; Sierpinska, 1992; Tall & Vinner, 1981; Vinner, 1983).

Such definition requires the knowledge of the set of real numbers, whose treatment is complex and rarely received by students. They are asked to operate on objects of which they normally have a partial understanding, using calculation algorithms not well mastered in their generality and accepted by "leap of faith" (think of the calculation of the square root of a number). The scarce operational content of numerical calculation emerges in all its gravity when they are requested to operate on the functions themselves or to imagine less standard functions, and this also occurs in university students, in which a greater knowledge is assumed; we therefore share the opinion of those who believe that the set-theoretic definition of function and the study of real functions of a real variable should not represent the starting point, but the result of a educational process aimed at building good mental images of the concept also through the reconstruction of its historical evolution.

To highlight the reasons why Descartes's algebraic approach provides a more immediate understanding, we may consider the construction of the $\sqrt{x}$ segment starting from the $x$ segment. This construction, from a geometric point of view, is that of a segment such that the square built on it is equal to the rectangle built on the starting segment $x$ and on a fixed segment, which we can call "unit".

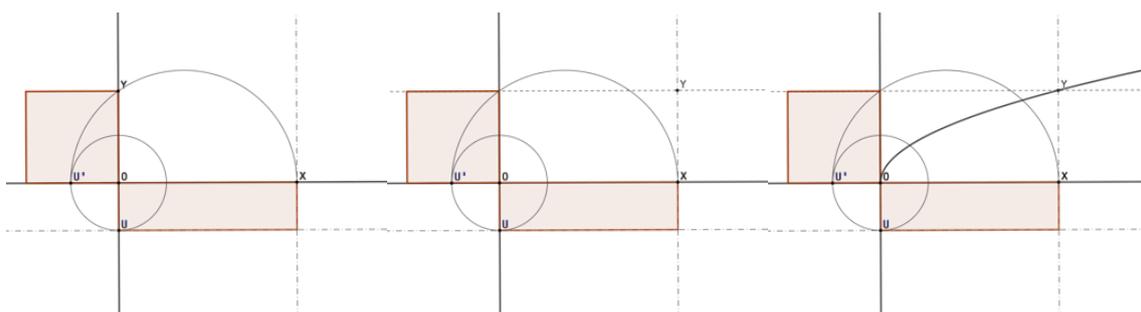

*Figure 2.* Given the rectangle $OUX$, to construct the equivalent square, consider the mean proportional $OY$ between the sides $OU$ and $OX$ of the rectangle (figure on the left). Let $XY$ be the segment parallel and congruent to $OY$ (figure in the center). As $X$ varies, the locus $Y$ of the segments $XY$ describes a curve, produced by a geometric covariation, expressed by the given construction. In Descartes's algebra, the covariation is described by the equation $y = \sqrt{x}$ where we have indicated with $x$ and $y$ the variable segments $OX$ and $OY$ respectively.

The construction of the root segment $OY$ is the same for all segments $OX$, and the variation of $OX$ occurs continuously, imagining to drag the point $X$ on the straight line $OX$, as physically can be done with a dynamic geometry software such as GeoGebra (Hoenwarter, 2002).

To carry out the same operation with numbers, by transforming the numerical variable $x$ into $y = \sqrt{x}$, the continuity of the operation is much less concrete, having to go through calculations of a very different nature, as for $\sqrt{1}, \sqrt{2}, \sqrt{\frac{\pi}{4}}$. Imagining a continuity behind these operations requires an effort of abstraction, obscuring the immediacy of the geometric variation (of $x$) and covariation (of $y$). Covariation is achieved with a geometric construction that presents the same difficulty for all segments and that can be expressed by an equation in the algebra of segments.

## 3. The algebra of segments

Descartes in his Géométrie begins by showing the constructions, with ruler and compass, of the sum, difference, product, quotient and square root, starting from any two segments. He is the first to understand the usefulness of thinking about the product and the quotient between two segments as a segment. In Euclid's Elements the product between two segments is a rectangle, while the quotient is a ratio, in the sense specified by Eudoxus and explained in Book V. The simple idea of basing constructions on a fixed segment (the unit segment) allows us to define a geometric algebra much more effective compared to that used by the Greeks, applying which it is possible to mathematically describe complicated covariation mechanisms without resorting to numbers.

We read the words that Descartes (1637) uses to introduce the fundamental operations in the first book of the Géométrie.

> *All the problems of Geometry can easily be reduced to such terms that a knowledge of the lengths of certain straight lines is sufficient for its construction. And just as arithmetic consists of only four or five operations, namely Addition, Subtraction, Multiplication, Division and Extraction of Root (which can be considered a sort of Division), so in geometry to find required lines it is merely necessary to add or subtract other lines; or else, taking one line, which I shall call unit in order to relate it as closely as possible to numbers, and which can be in general chosen arbitrarily, and having given two other lines, it is possible to find a fourth line which shall be to one of these two as the other is to unit (which is the same as multiplication); or, again, to find a fourth line that shall be to one of the given lines as unit is to the other (which is the same thing as division); or finally, find one, two, or several mean proportional between the unit and some other lines (which is the same as extracting the square root, cube root, etc. of the given line). And I shall not hesitate to introduce these arithmetical terms into geometry, for the sake of a greater clearness.*
> *Multiplication. For example, let AB be taken as unit and let it be required to multiply BD by BC; I have only to join the points A and C, and draw DE parallel to CA, then BE is the product of this Multiplication.*
> *Division. Or if it is required to divide BE by BD, I join the points B and E, and draw AC parallel to DE; then BC is he result of the division.*
> *Extraction of the square root. Or if the square root of GH is required I add, along the same straight line, FG equal to unit; then, bisecting FH at K, I describe the circumference FIH about K as a centre, and draw from G a perpendicular and extend it to I, and GI is the required root I do not speak here of cube root, or other roots, since I shall speak more conveniently of them later* (Descartes, 1632, pp. 297–298).

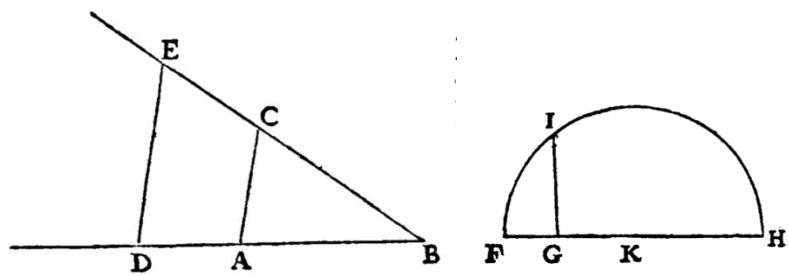

*Figure 3.* Thales's theorem allows us to interpret the product and ratio of two segments as a segment (figure on the left (Descartes, 1637, p. 298)). The mean proportional between a segment and the unit represents the root of the segment (figure on the right, (Descartes, 1637, p. 298)).

Combining these constructions, starting from a segment $x$ we can construct a segment $y$ with ruler and compass whose dependence on $x$ is expressed inalgebra of segments through a rational expression $y = f(x)$ (function of $x$ or square root of $x$). We can transport $y$ in order to make an extreme coincide with the varying extreme of $x$, perpendicular to it. The other extreme draws a curve that can be very effectively visualized with GeoGebra: the graph of the function $f$. In our opinion, the constructions with ruler and compass link the graph to the expression in a more natural way than arithmetic operations on numbers.

Using algebra of segments, Descartes explains how a geometric construction corresponds to the solution of a system of polynomial equations. If the system is indeterminated, the locus is obtained as the segment varies. Our intent is to show how it is possible, through a series of elementary geometric constructions that can be encoded in algebraic expressions, to specify the connection between two variable segments. We believe that, from the student's point of view, the use of segments instead of numbers and ruler and compass constructions instead of numerical operations gives a more concrete meaning to the continuous variation and covariation of mathematical objects; we hope this can provide a better understanding of the concept of real function of real variable and of general functions.

## 4. A proposal for an educational path

According to the previous reflections, we propose an educational path to be carried out in a secondary school class. In this work we limit to suggest the main stages without going into the implementation details that will be exhaustively treated elsewhere. The proposal has a laboratory character and uses the GeoGebra dynamic geometry software (Hoenwarter, 2002).

### 4.1 The graph of the function "square" $y = x^2$

We consider a variable $OX$ segment and on it we construct a square, as indicated in Euclid I. 46. By varying the point $X$, the square produced by the construction adapts continuously to the variation of the segment, as it's easily observed with a dynamic geometry software.

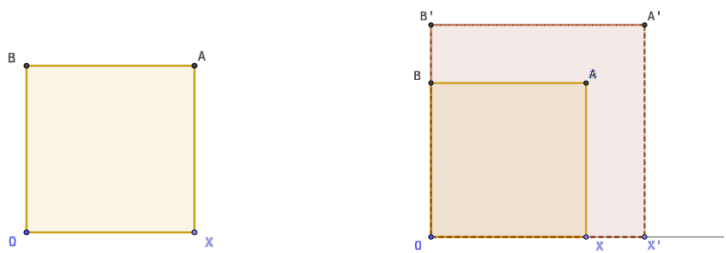

*Figure 4.* In the figure on the left, the construction of a square. In the figure on the right, the largest square is obtained by dragging the construction with the extreme variable $X$ of the $OX$ segment.

In the variable square, all four sides vary simultaneously. Let us think of fixing a point $U$ on the half-line $OX$ and then building on $OU$ the rectangle equivalent to the square with side $OX$. To this end, we trace from $X$ the parallel to $BU$ which intersects the line $OB$ in $C$. The theorem of Thales tells us that $OX: OU = OC: OB$, or rather the rectangle sought is that of sides $OU$ and $OC$.

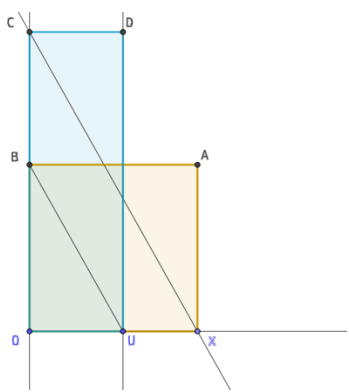

*Figure 5.* Use of Thales's theorem to construct a rectangle equivalent to a square.

The choice of $U$ therefore allowed us to represent the variation of the extension of the square on $OX$ with a rectangle that has a fixed and a variable segment. Hence we have concentrated the variation of the square on $OX$ in that of a single segment $OU$, which we can carry on $XY$, perpendicular to the half-line $OU$.

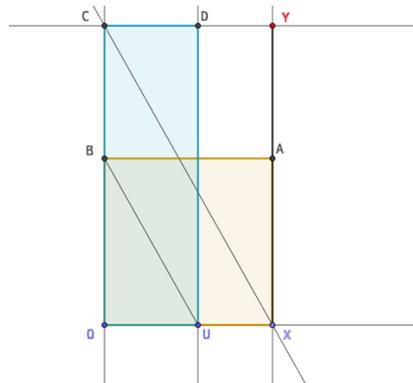

*Figure 6.* Transport of one side of a rectangle having the other fixed side perpendicular to a variable segment.

The locus of points $Y$ is therefore able to represent the variation of the square $OX$.

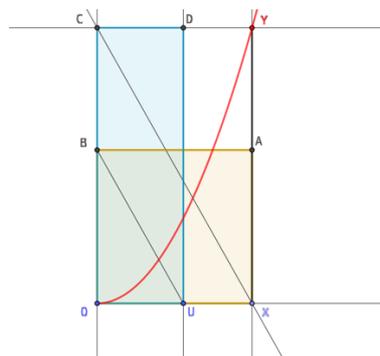

*Figure 7.* Locus generated by the extremes of the $XY$ segment by varying the extreme $X$ of the variable segment $OX$.

The proportion $y:x = x:1$ (where $1 = OU$, $x = OB$ and $y = OY = OY'$), which describes the "symptom" of the locus of points, can be expressed with the equation $y = x^2$.
    We begin with the geometric construction of the square of a segment $OX$ in which $O$ is fixed and $X$ is free to vary, and we show how in Descartes's algebra the square of a segment is itself a segment. We introduce the idea of how to proceed to visualize a covariation between variable segments: starting from the half-line of origin $O$, on which the unity $OU$ is fixed, we construct the variable segment $O$. Then the $OY$ segment constructed starting from the $OX$ segment is transported on the perpendicular line passing through $X$, in order to "separate" the variations and thus generate the locus (curve) described by the point $Y'$. The "symptom" of the curve constructed with this geometric covariation can be expressed with the proportion $y:x = x:1$, (where $1 = OU$, $x = OX$ e $y = OY = OY'$) from which we get the equation $y = x^2$. Thanks to the introduction of the unit segment, each relationship and each product can be represented with a segment. The equations lose the homogeneity of the proportions but are more flexible.

## 4.2 The graph of the function "cube" $y = x^3$

In Descartes' algebra we can also represent the cube of an $OX$ segment with a segment, without going through a three-dimensional construction of the cube. Multiplying the segment $OY = x^2$ obtained in the previous construction by the segment $x$ and again using the Thales's theorem, we can construct the segment $y = x^3$ whose locus, varying $X$ and using the language of functions, describes the graph of $y = x^3$.

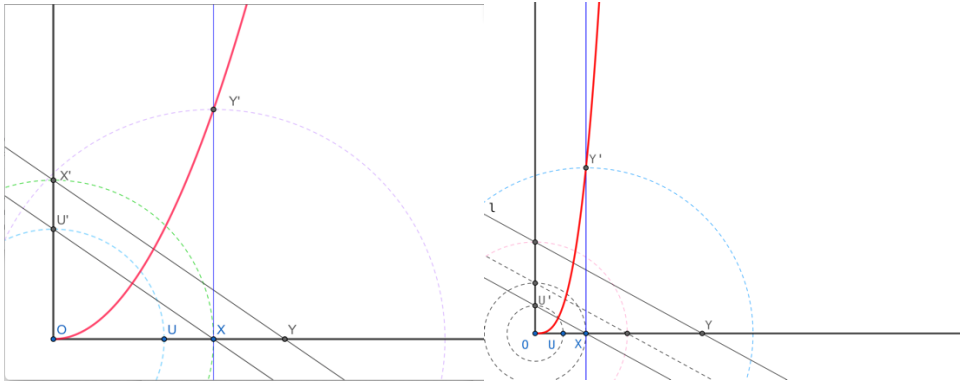

*Figure 8.* In the figure on the left, the construction of the segment $XY'$ which represents $OX^2$ in Descartes's algebra, is carried out in a slightly different way from the one we have described in subsection 4.1, which is better adapted to iteration. We consider the $OU$ on the orthogonal semi-axis, then we intersect the circumference centered in $O$ and passing through $U$ with the semi-axis at the point $U'$. Similarly, the $OX$ segment is shown on $OX'$. By Thales's theorem, the parallel for $X'$ to $XU'$ cuts the semi-axis $OX$ at $Y$ such that $OY = OX^2$. The extreme $Y'$ of the segment $XY'$ congruent to $OY$ on the perpendicular to $OX$ in the semi-plane containing $X'$ (which can be constructed with ruler and compass using Euclid I.2) describes the graph of $y = x^2$. In the figure on the right, to multiply the segment $x^2$ obtained in the previous step, it is sufficient to bring its end to the $OU'$ semi-axis with the compass and trace the parallel to $U'X$ from the point obtained. Its intersection with the $OX$ semi-axis determines the segment $OY = OX^3$. The extreme $Y'$ of the segment $XY'$ congruent to $OY$ describes the graph of $y = x^3$.

## 4.3 the graph of the function "hypercube" $y = x^4$

In Descartes's algebra we can easily construct the graph of the function $y = x^4$, which geometrically is not associated with any elementary figure in the same sense in which the square is associated with $y = x^2$ and the cube with $y = x^3$. The segment $x^4$ can be obtained either by multiplying $x^3$ by $x$, or by multiplying $x^2$ by $x^2$.

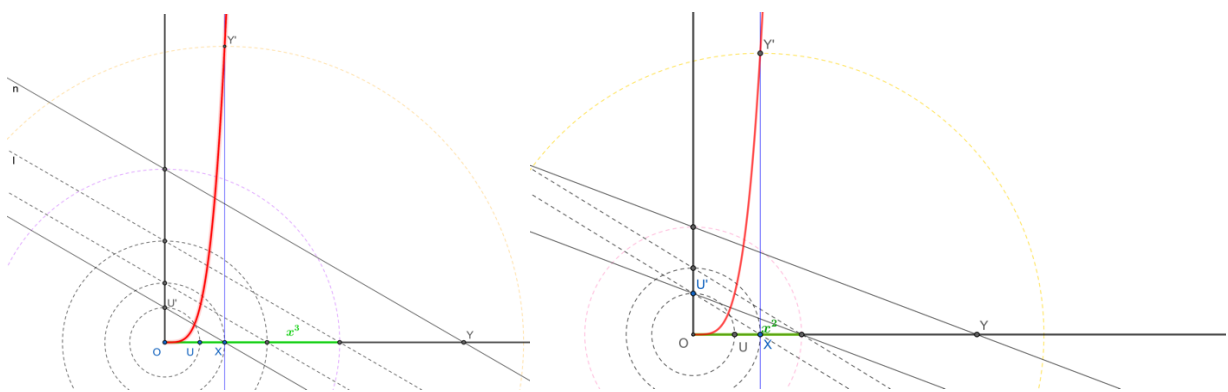

*Figure 9.* In the left figure, the construction of the segment $y = x \cdot x^3$ In the figure on the right, the construction of the segment $y = x^2 \cdot x^2$.

### 4.4 The graph of a polynomial function $y = a_0 + a_1x + a_2x^2 + \ldots + a_nx^n$

After showing the procedure for constructing the powers of a segment $x$, it should be clear how it is possible to construct with ruler and compass the segment represented by the combination $y = a_0 + a_1 + a_2x^2 + \ldots + a_nx^n$.

The coefficients $a_i$ of the combination are also represented in Descartes's algebra by segments, which can be drawn in such a way as to be able to geometrically interpret the dependence on the parameters of the graph of the polynomial function simply dragging the ends of these segments.

### 4.5 The graph of the function "reciprocal" $y = \frac{1}{x}$

Also the construction $y = \frac{1}{x}$ (and so the graph of the function) can be easily made:

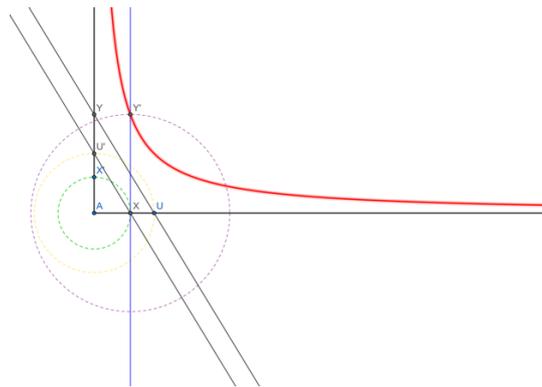

*Figura 10.* Construction of the segment $y = \frac{1}{x} = XY'$ starting from the segment $x = AX$, using the Thales's theorem.

### 4.6 The graph of the function "square root" $y = \sqrt{x}$

As mentioned in paragraph 3, Descartes also teaches to construct the square root of a segment $x$, and therefore the graph of the function $y = \sqrt{x}$.

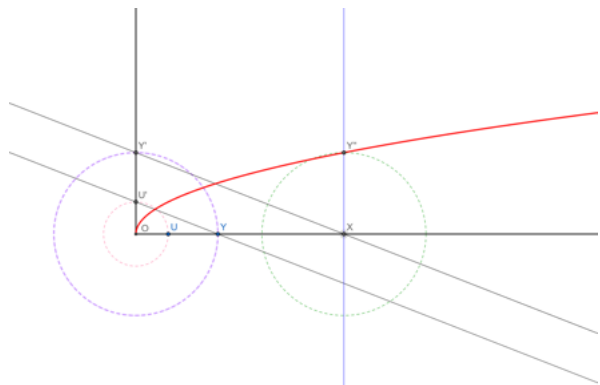

*Figure 11.* Construction of the segment $y = \sqrt{x} = XY''$ starting from the segment $x = OY$.

### 4.7 Conclusion

Putting together what we have seen so far, it is not difficult to imagine how to construct the graph of a rational function of the type $y = f(x)$ or $y = f(\sqrt{x})$. We need to keep in mind that the loci we have built are not constructed with ruler and compass. Such constructions produce a single point of the locus, but not the continuous variation of that point, i.e. the whole curve.

With the "locus command" of GeoGebra it is also possible to draw the graph of the inverse function of a function whose locus has already been constructed. Let's consider a point $Y$ on the graph of a function $y = f(x)$ and construct its symmetric, with ruler and compass, with respect to the bisector of the two semi-axis. The locus of these points draws the graph of the inverse function.

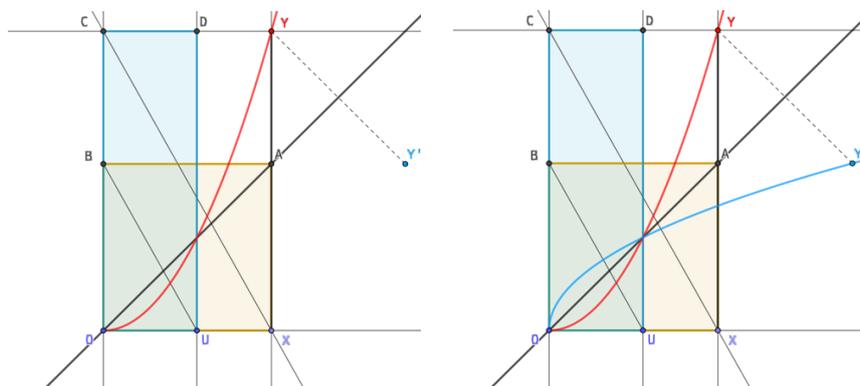

*Figura 12*. The locus of symmetrical points (with respect to the bisector of the two semi-axis) of the points of the graph of $y = f(x)$ is the graph of the function $y = f^{-1}(x)$.

However, has to be noted that this construction does not create a ruler and compass construction of a segment $z = f^{-1}(y)$ starting from a segment $y$, but only starting from a segment $y = f(x)$. The construction of $f(y)$, given $y$, is not generally possible. For example, it is possible to construct, with ruler and compass, any finite number of points on the graph of the cube root but not the cube root of a segment.

The above brings out an aspect of the graph of a function as a mathematical object, on which it is possible to operate with constructions of a higher level than those that have allowed us to "construct the points of the graph of a function": a fundamental step in the "reification " (Sfard, 1992) of a concept. Finally, we observe how we have always limited the construction of the graph to the quadrant between two orthogonal semi-axis having a common endpoint. Descartes's algebra is not an algebra with oriented segments, and therefore the geometric treatment of negative quantities is excluded. We hope to come back to this point in a later work.